\documentclass[12pt,reqno,a4paper]{amsart}
\usepackage{verbatim,rcs,cite}
\usepackage{amssymb,amsfonts,latexsym,mathrsfs}

\RCS $Revision:  $ \RCSdate $Date: 2008/01/29 13:25:27 $

\title[The Minimal Degree for a Class of Groups]{The Minimal Degree for a Class of Finite Complex Reflection Groups}

\date{}

\author{Neil Saunders}

\dedicatory{\upshape
School of Mathematics and Statistics \\
University of Sydney, NSW 2006, Australia\\[.5em]
{\it E-mail address:} \ \texttt{neils@maths.usyd.edu.au}
}

\usepackage{amsmath}
\usepackage{amssymb}
\usepackage{amsthm}
\usepackage{mathrsfs}
\usepackage[all]{xy}
\usepackage{pstricks}
\newtheorem{theorem}{Theorem}[section]
\newtheorem{lemma}[theorem]{Lemma}
\newtheorem{proposition}[theorem]{Proposition}
\newtheorem{corollary}[theorem]{Corollary}

\theoremstyle{remark}

\theoremstyle{definition}
\newtheorem{definition}[theorem]{Definition}

\DeclareMathOperator{\core}{core} \DeclareMathOperator{\soc}{soc} \DeclareMathOperator{\codim}{codim}

\begin{document}



\maketitle

\begin{abstract}
We calculate the minimal degree for a class of finite complex reflection groups $G(p,p,q)$, for $p$ and $q$ primes and establish relationships between minimal degrees when these
groups are taken in a direct product.
\end{abstract}

\section{Introduction}
The minimal faithful permutation degree $\mu(G)$ of a finite group $G$ is the least non-negative integer $n$ such that $G$ embeds in the
symmetric group $Sym(n)$. It is well known that $\mu(G)$ is the smallest value of $\sum_{i=1}^{n}|G:G_i|$ for a collection of subgroups
$\{G_1,\ldots, G_n \}$ satisfying $\bigcap_{i=1}^{n}\core(G_i)=\{1\}$, where $\core(G_i)=\bigcap_{g \in G}G_{i}^{g}$. \vspace{12pt}

We will often denote such a collection of subgroups by $\mathscr{R}$ and refer it as the representation of $G$. The elements of $\mathscr{R}$
are called {\it transitive constituents} and if $\mathscr{R}$ consists of just one subgroup $G_0$ say, then we say that $\mathscr{R}$ is
transitive and that $G_0$ is core-free. \vspace{12pt}

The study of this are dates back to Johnson \cite{J71} where he proved that one can construct a minimal faithful representation $\{G_1,\ldots,
G_n \}$ consisting entirely of so called {\it primitive} groups. These are groups which cannot be expressed as the intersection of groups that
properly contain it. \vspace{12pt}

Here we give a theorem due to Karpilovsky \cite{K70}, which also serves as an introductory example. We will make use of this theorem later and
the proof of it can be found in \cite{J71} or \cite{S05}.
\begin{theorem} \label{theorem:abelian}
Let A be a finite abelian group and let $A \cong A_1 \times \ldots \times A_n$ be its direct product decomposition into non-trivial cyclic
groups of prime power order. Then $$\mu(A)= a_1 + \ldots + a_n,$$ where $|A_i|=a_i$ for each $i$.
\end{theorem}

The aim of this paper is calculate $\mu(G(p,p,q))$ where $p$ and $q$ are prime numbers and $G(p,p,q)$ is the monomial complex reflection group
(see Section \ref{section:calculations} for the definition). \vspace{12pt}

When $q=2$, the group $G(m,m,2)$ is the dihedral group of order $2m$ ($m$ an integer) and in \cite[Proposition 2.8]{EP88}, Easdown and Praeger
calculated the minimal degree for the dihedral groups. Specifically, they proved the following:

\begin{proposition}
For any integer $k=\prod_{i=1}^{m}p_{i}^{\alpha_i} >1$, with the $p_i$ distinct primes, define $\psi(k)=\sum_{i=1}^{m}p_{i}^{\alpha_i}$, with
$\psi(1)=0$. Then for the dihedral group $D_{2^{r}n}$ with $n$ odd, we have
\begin{displaymath}
\mu(D_{2^{r}n})= \left\{ \begin{array}{lll}
2^r &\text{if} \quad n=1, 1\leq r \leq 2 \\
2^{r-1} & \text{if} \quad  n=1, r>2 \\
\psi(n) & \text{if} \quad  n>1, r=1 \\
2^{r-1} + \psi(n) & \text{if} \quad  n>1, r>1.\\
\end{array} \right.
\end{displaymath}
\end{proposition}

When $p=2$ and $q=3$, the group $G(2,2,3)$ is isomorphic to $Sym(4)$ and so $\mu(G(2,2,3))=4$. So for the rest of this article, we will assume
that $p$ and $q$ are odd primes. The main result of this article (see Theorem \ref{theorem:result1}, Theorem \ref{theorem:result2} and
Proposition \ref{proposition:result3}) is that

\begin{displaymath}
\mu(G(p,p,q))= \left\{ \begin{array}{lll}
pq &\text{if} \quad p<q \quad \text{or,} \quad p>q \geq 5 \quad \text{or,} \quad q=3 \ \text{and} \ p \equiv 2 \ \text{mod} \ 3\\
p^2 &\text{if} \quad p=q\\
2q & \text{if} \quad  q=3 \ \text{and} \ p \equiv 1 \ \text{mod} \ 3.\\
\end{array} \right.
\end{displaymath}

In Section \ref{section:directproduct} we continue the theme of Johnson \cite{J71} and Wright \cite{W75}, of investigating the inequality
\begin{equation} \mu(G \times H) \leq \mu(G) + \mu(H) \label{eq:directsum} \end{equation} for finite groups $G$ and $H$. Johnson and Wright
first investigated under what conditions equality holds in (\ref{eq:directsum}), (see Section \ref{section:directproduct} for more details).
However examples for when the inequality is strict are not well known. \vspace{12pt}

We show that standard wreath $C_p \wr Sym(q)$ is the internal direct product of $G(p,p,q)$ and its non-trivial centralizer in $C_p \wr Sym(q)$
and furthermore,
$$\mu(C_p \wr Sym(q))=\mu(G(p,p,q))=pq,$$ so we get a family of examples of a strict inequality in (\ref{eq:directsum}).

\section{Background}
\subsection{Definitions Relating to the Socle} The socle of a finite group $G$ is defined to be the subgroup
generated by all of the minimal normal subgroups of $G$. It follows that this is the direct product of all of the minimal normal subgroups of
$G$. That is,
$$\soc(G)=\prod_{i}^{d}M_i,$$ where each $M_i$ is a minimal normal subgroup of $G$ and $d$ is the number of minimal normal subgroups.

\begin{definition}
For a finite group $G$ define $\dim(G)$ to be the number of minimal normal subgroups of $G$. For a subgroup $H$ of $G$ define $\dim_{G}(H)$ to
be the number of minimal normal subgroups of $G$ that are contained in $H$. Also define $\codim_{G}(H)=\dim(G) - \dim_{G}(H)$. By convention,
the trivial group is defined to have dimension $0$.
\end{definition}

\subsection{A Note on Cyclotomic Polynomials and Roots of Unity}

We now state a result from \cite[Theorem 10.2.4]{FT95} which will be needed later.

\begin{theorem}
Let $p$ and $q$ be distinct primes and let $S$ be a splitting field for $x^q-1$ over the field $\mathbb{F}_p$. Then $S=\mathbb{F}_{p^{o}}$ where
$o$ is the multiplicative order of $p$ mod $q$.
\end{theorem}

Given this theorem, we immediately have the following.

\begin{corollary} \label{corollary:splittingfield}
Let $p$ and $q$ be distinct primes. Then $\mathbb{F}_p$ is the splitting field for $x^q-1$ if and only if $p \equiv 1$ mod $q$.
\end{corollary}

This corollary shows that a $q$-th root of unity exists in $\mathbb{F}_p$ if and only if $p \equiv 1$ mod $q$. \vspace{12pt}

The next definition and result is taken from \cite{LN97}.

\begin{definition}
For $r$ a prime number, the polynomial $$Q_r(x)=1+x+x^2+\ldots+x^{r-1}$$ is called the $r$-th {\it cyclotomic polynomial}. The roots of this
polynomial are non-trivial $r$-th roots of unity.
\end{definition}

\begin{theorem} \label{theorem:irreduciblesubmodules}
Let $\mathbb{F}_q$ be a finite field of $q$ elements and let $n$ be a positive integer coprime to $q$. Then the polynomial $Q_{n}(x)$ factors
into $\frac{\phi(n)}{d}$ distinct monic irreducible polynomials in $\mathbb{F}_q[x]$ of the same degree $d$ where $d$ is the least positive
integer such that $q^d \equiv 1$ mod $n$.
\end{theorem}

Thus for $r$ a prime, $Q_r(x)$ splits into $\frac{r-1}{d}$ monic irreducible factors where $d$ is the multiplicative order of $r$ in the group
of units $(\mathbb{Z}/n\mathbb{Z})^*$.

\subsection{The Groups $G(m,p,n)$} In this section we follow the notation of \cite{TL07}. \vspace{12pt}

Let $m$ and $n$ be positive integers, let $C_m$ be the cyclic group of order $m$ and $B=C_m \times \ldots \times C_m$ be the direct product of
$n$ copies of $C_m$. For each divisor $p$ of $m$ define the group $A(m,p,n)$ by
$$A(m,p,n)= \{(\theta_1,\theta_2,\ldots, \theta_n) \in B  \ | \ (\theta_1\theta_2\ldots\theta_n)^{m/p}=1 \}.$$
It follows that $A(m,p,n)$ is a subgroup of index $p$ in $B$ and the symmetric group $Sym(n)$ acts naturally on $A(m,p,n)$ by permuting the
coordinates. \vspace{12pt}

$G(m,p,n)$ is defined to be the semidirect product of  $A(m,p,n)$ by $Sym(n)$. It follows that $G(m,p,n)$ is a normal subgroup of index $p$ in
the wreath product $C_m \wr Sym(n)$ and thus has order $m^nn!/p$. \vspace{12pt}

It is well known that these groups can be realized as finite subgroups of $GL_n(\mathbb{C})$, specifically as $n \times n$ matrices with exactly
one non-zero entry, which is a complex $m$-th root of unity, in each row and column such that the product of the entries is a complex $(m/p)$th
root of unity. Thus the groups $G(m,p,n)$ are sometimes referred to as monomial reflection groups. For more details on the groups $G(m,p,n)$,
see \cite{OT92}, \cite{C76}.

\section{Calculation of Minimal Degrees} \label{section:calculations}

Throughout this section, $p$ and $q$ will be distinct odd primes. Recall that $G(p,p,q)=A(p,p,q) \rtimes Sym(q)$ where
$$A(p,p,q)=\{(\theta_1,\theta_2,\ldots, \theta_q) \in B \ | \ \theta_1\theta_2\ldots\theta_q=1 \},$$ which is isomorphic to the direct product
of $q-1$ copies of the cyclic group of order $p$, $C_p$. Hence
$$G(p,p,q) \cong \underbrace{(C_p \times \ldots \times C_p )}_{q-1} \rtimes Sym(q).$$ From now on, we will let $G$ denote the group $G(p,p,q)$
and $A$ denote the group $A(p,p,q)$. Since $A$ is a proper subgroup of $G$ we have by Theorem \ref{theorem:abelian}, $\mu(A)=p(q-1)$ and so
$p(q-1) \leq \mu(G)$. On the other hand, $G$ is a proper subgroup of the wreath product $C_p \wr Sym(q)$ which can easily be seen to have
minimal degree $pq$, thus
\begin{equation} p(q-1) \leq \mu(G) \leq pq. \end{equation} \label{eq:inequality}

Since $p$ is prime, we may treat $B$ as a $Sym(q)$-module over the finite field $\mathbb{F}_p$ with basis $\theta_1, \ldots, \theta_q$. In this
setting, $A$ is the codimension one submodule of $B$ consisting of the elements
$$\prod_{i=1}^{q}\theta_{i}^{\lambda_i} \quad \text{with} \quad \sum_{i=1}^{q}\lambda_{i}=0.$$ It is clear that
$c_1=\theta_1\theta_2^{-1}, \ldots, c_{q-1}=\theta_{q-1}\theta_{q}^{-1}$ is a basis for $A$. Let $b=(1 \ 2 \ldots \ q)$ be the $q$-cycle in
$Sym(q)$. Then $b$ acts on the basis for $A$ in the following manner
\begin{equation} c_{1}^b=c_2, \ c_{2}^b=c_3, \ldots, c_{q-2}^b=c_{q-1}, \ c_{q-1}^b=(c_1c_2\ldots c_{q-1})^{-1}.
\label{equation:baction} \end{equation} The matrix for the action of $b$ with respect to this basis is the companion matrix

\begin{displaymath}
\left( \begin{array}{cccccc}
0       &   1       &   0       &   \ldots  &   0      &   0 \\
0       &   0       &   1       &   \ldots  &   0      &   0 \\
0       &   0       &   \ddots  &   \ddots  &   \vdots &    \vdots\\
\vdots  &   \vdots  &   \ddots  &           &   0      &   1  \\
-1      &  -1       &   -1      &   \ldots  &   -1     &    -1 \\
\end{array} \right)
\end{displaymath}

The minimal polynomial for this action is the cyclotomic polynomial $Q_{q}(\lambda)=1+\lambda+\lambda^2+\ldots+\lambda^{q-1}$. We can make an
important observation at this point.

\begin{proposition} \label{proposition:submodule}
$b$ does not commute with any non-trivial element of $A$.
\end{proposition}
\begin{proof}
Suppose there is an element of $A$ which commutes with $b$. Then this element corresponds to an eigenvector for $b$ with eigenvalue $1$.
However, this implies that $1$ is a solution to the characteristic polynomial $Q_q(\lambda)$, a contradiction.
\end{proof}

We will make greater use of the action of $b$ on $A$ and the cyclotomic polynomial $Q_q(\lambda)$ in later sections.\vspace{12pt}

If $e_1, \ldots, e_q$ is the dual basis for $B$ where $e_i(\theta_j)=\delta_{ij}$, then the basis of $A^*$ dual to $\{c_i\}$ is $f_1=e_1,
f_2=e_1+e_2, \ldots, f_{q-1}=e_1+\ldots+e_{q-1}$. It easily follows that $f_i(c_j)=\delta_{ij}$. Observe that the functional $e_1+\ldots+e_q=0$
when restricted to $A$, since $A$ is the subspace of $B$ whose elements have the property that the sum of their coefficients is zero. We write
$A^*$ and $B^*$ additively and we have $A^*=B^*/\mathbb{F}_p(e_1+\ldots+e_q)$.\vspace{12pt}

\subsection{$\mu(G(p,p,q))$ for $p > q$}\label{section:p>q}

In this subsection we assume that $p > q$ and exploit the action of $Sym(q)$ on $A$ to prove that every minimal faithful representation of $G$
is given by a core-free subgroup.\vspace{12pt}

The following is a well known result from modular representation theory.

\begin{proposition}
$Sym(q)$ acts irreducibly and faithfully on $A$.
\end{proposition}
\begin{proof}
We show that the submodule generated by an arbitrary non-trivial element is the whole of $A$. Let $w=\prod_{i=1}^{q}\theta_{i}^{\lambda_i}$ be a
non-trivial element of $A$ so that $\sum_{i=1}^{q}\lambda_{i}=0$. It is enough to prove that we can obtain the basis elements
$c_1=\theta_1\theta_2^{-1},\ldots, c_{q-1}=\theta_{q-1}\theta_{q}^{-1}$ of $A$ via the action of $Sym(q)$ on $w$.\vspace{12pt}

Fix a non-zero $\lambda_i$. There is another non-zero $\lambda_j$ such that $\lambda_i-\lambda_j \neq 0$. For if $\lambda_i-\lambda_k=0$ for all
$k$, then we would have $\lambda_i=\lambda_k$ for all non-zero $\lambda_k$ and so
$$w=(\prod_{i \in I}\theta_i)^{\lambda_i}$$ where $I$ is a subset of $\{1,2,\ldots,q\}$. So $\sum_{i \in
I}\lambda_i=|I|\lambda_i=0$ in $\mathbb{F}_p$. However since $p>q$, this implies that $\lambda_i=0$ contradicting that $w$ is
non-trivial.\vspace{12pt}

Choose two non-zero $\lambda_i$ and $\lambda_j$ with $\lambda_i-\lambda_j \neq 0$. Then applying the transposition $(i \ j)$ to $w$ we have
$$w^{(i \ j)}=\theta_{1}^{\lambda_1}\ldots\theta_{i}^{\lambda_j}\ldots\theta_{j}^{\lambda_i}\ldots\theta_{q}^{\lambda_q},$$ so $$w(w^{(i \
j)})^{-1}=\theta_i^{\lambda_i-\lambda_j}\theta_j^{\lambda_j-\lambda_i}=(\theta_i\theta_j^{-1})^{\lambda_i-\lambda_j}.$$ Therefore,
$\theta_i\theta_j^{-1}$ is contained in $A$ and by applying the appropriate permutation to it, we can obtain all the basis elements $c_1,
\ldots,c_{q-1}$ as required. So $Sym(q)$ acts irreducibly on $A$.\vspace{12pt}

Now suppose that the action of $Sym(q)$ on $A$ has a kernel. This kernel must be normal subgroup of $Sym(q)$ and since $q \neq 4$, the only
possibility is the alternating group $Alt(q)$. However, the $q$-cycle $b$ is an even permutation which does not commute with any non-trivial
element of $A$. Therefore $Sym(q)$ acts faithfully on $A$.
\end{proof}

\begin{corollary} \label{corollary:minnorsubgroup}
$A$ is the unique minimal normal subgroup of $G$.
\end{corollary}
\begin{proof}
Certainly $A$ is a normal subgroup of $G$ and since $Sym(q)$ acts irreducibly on it, it is a minimal normal subgroup.\vspace{12pt}

Suppose $N$ is a non-trivial normal subgroup of $G$ which does not contain $A$. By minimality of $A$ we must have $A \cap N=\{1\}$. It follows
that $AN=A \times N$, that is $A$ and $N$ pairwise commute. \vspace{12pt}

Let $a \in A$ and $n=a'\sigma \in N$, where $a' \in A$ and $\sigma \in Sym(q)$. Then $n=a^{-1}na$ so $a'\sigma=a^{-1}a'\sigma a= a'a^{-1}\sigma
a$, and so $\sigma=a^{-1}\sigma a$. That is, $\sigma$ commutes with $a$ and so $\sigma$ is contained in the kernel of the action of $Sym(q)$ on
$A$. Therefore $\sigma$ is trivial and $n=a'$ contradicting that $A \cap N=\{1\}$.\vspace{12pt}

Therefore $A$ is contained in every non-trivial normal subgroup of $G$ and is thus the unique minimal normal subgroup of $G$.
\end{proof}

It follows now that any minimal faithful representation of $G$ must be transitive, that is, given by a single core-free subgroup.\vspace{12pt}

Let $L$ be a core-free subgroup of $G$ such that $|G:L|=\mu(G)$. Recall we have the upper and lower bounds $p(q-1) \leq \mu(G) \leq pq$. Since
$|G|=p^{q-1}q!$ the upper and lower bounds above simplify to $$q-1 \leq \frac{p^{q-2}q!}{|L|} \leq q.$$ Observe that $|L|$ is not divisible by
$p^{q-1}$ since $L$ is core-free, so $\frac{p^{q-2}q!}{|L|}$ is indeed an integer. Therefore if $\mu(G) \neq pq$ then $\mu(G)=p(q-1)$. We aim to
show that this is not the case (except when $q=3$ and $p \equiv 1$ mod $3$) via an argument by contradiction. First we suppose that the odd
prime $q$ is at least $5$. \vspace{12pt}

Assume $\mu(G)=p(q-1)$ so that $|L|=p^{q-2}q(q-2)!$. Then since $p>q$ and $A$ is the unique Sylow $p$- subgroup of $G$, we have that $L$
contains a unique Sylow $p$-subgroup $A \cap L$ of order $p^{q-2}$. Moreover by Sylow's Theorem and since we are assuming $q$ is at least $5$,
$L$ contains an element of order $q$, which without loss of generality, we may assume this element is the $q$-cycle $b=(1 \ 2 \ldots q)$ and $L$
also contains and element of order $3$, $x$ say. \vspace{12pt}

Now it is cleat that we may simultaneously treat $A$ as a $\langle b \rangle$-module and as an $\langle x \rangle$-module, and in this respect,
$A \cap L$ as a codimension-$1$ submodule of $A$ for both these cyclic groups $\langle b \rangle$ and $\langle x \rangle$. Therefore we can
consider $A \cap L$ as the kernel of a linear functional $f$ on $A$. Upon expressing $f$ in terms of the dual basis for $A^*$ we have $f=a_1f_1+
\ldots +a_{q-1}f_{q-1}$ for some $a_i \in \mathbb{F}_p$ not all zero. Therefore we may write $$A \cap L = \{c_{1}^{\alpha_1}c_2^{\alpha_2}\ldots
c_{q-1}^{\alpha_{q-1}} \ | \ a_1\alpha_1 + a_2\alpha_2 + \ldots + a_{q-1}\alpha_{q-1}=0 \}.$$

Since $b$ preserves $A \cap L$, the condition $c_{1}^{\alpha_1}c_2^{\alpha_2}\ldots c_{q-1}^{\alpha_{q-1}} \in A \cap L$ must be equivalent to
$(c_{1}^{\alpha_1}c_2^{\alpha_2}\ldots c_{q-1}^{\alpha_{q-1}})^b \in A \cap L$. Now
\begin{eqnarray*}
(c_{1}^{\alpha_1}c_2^{\alpha_2}\ldots c_{q-1}^{\alpha_{q-1}})^b &=& c_{2}^{\alpha_1}c_3^{\alpha_2}\ldots c_{q-1}^{\alpha_{q-2}}(c_1c_2\ldots c_{q-1})^{-\alpha_{q-1}} \\
 &=& c_{1}^{-\alpha_{q-1}}c_2^{\alpha_1-\alpha_{q-1}}\ldots c_{q-1}^{\alpha_{q-2}-\alpha_{q-1}}.
 \end{eqnarray*}
So the condition for $(c_{1}^{\alpha_1}c_2^{\alpha_2}\ldots c_{q-1}^{\alpha_{q-1}})^b$ to be in $A \cap L$ is:
$$a_1(-\alpha_{q-1})+a_2(\alpha_1-\alpha_{q-1})+ \ldots +a_{q-1}(\alpha_{q-2}-\alpha_{q-1})=0.$$ Rewriting this equation in terms of the
$\alpha_i$'s we have
$$a_2\alpha_1+a_3\alpha_2+\ldots+a_{q-1}\alpha_{q-2}+(-a_1-a_2-\ldots-a_{q-1})\alpha_{q-1}=0.$$ Upon equating the coefficients up to scalar
equivalence, we get
\begin{eqnarray*}
a_1 &=& \lambda a_2 \\
a_2 &=& \lambda a_3 \\
    &\vdots&            \\
a_{q-2} &=& \lambda a_{q-1} \\
a_{q-1} &=& \lambda(-a_1-a_2-\ldots-a_{q-1}) a_{q-1},
\end{eqnarray*}
for some non-zero $\lambda \in \mathbb{F}_p$. \vspace{12pt}

Hence we may express each $a_i$ as $\lambda^{q-1-i}a_{q-1}$ and so the last equation can be expressed as
$$(\lambda^{q-1}+\lambda^{q-2}+\ldots+1)a_{q-1}=0.$$ Now $a_{q-1} \neq 0$ since otherwise all of the $a_i$ would be zero, so it follows that $$\lambda^{q-1}+\lambda^{q-2}+\ldots+1=0.$$
Therefore, since $q$ is prime, $\lambda$ is a primitive $q$-th root of unity $\zeta_q$ which exists if and only if $p \equiv 1$ mod $q$ by
Corollary \ref{corollary:splittingfield}. So if $p$ is not congruent to $1$ modulo $q$, then we have our desired contradiction proving that
$\mu(G)=pq$. Therefore from now on, we assume $p \equiv 1$ mod $q$. \vspace{12pt}

The condition for the element $c_{1}^{\alpha_1}c_2^{\alpha_2}\ldots c_{q-1}^{\alpha_{q-1}} \in A \cap L$ to satisfy is now
$$\zeta_q^{q-2}\alpha_1+\zeta_q^{q-3}\alpha_2+\ldots+\zeta_q\alpha_{q-2}+\alpha_{q-1}=0,$$ and so
$f=\zeta_q^{q-2}f_1+\zeta_q^{q-3}f_2+\ldots+\zeta_qf_{q-2}+f_{q-1}$. Upon expressing each of the $f_i$'s in terms of the $e_i$'s we have
\begin{eqnarray*}
f &=& \zeta_q^{q-2}e_1+\zeta_q^{q-3}(e_1+e_2)+\ldots+\zeta_q(e_1+\ldots +e_{q-2})+(e_1+\dots+e_{q-1}) \\
  &=& (\zeta_q^{q-2}+ \zeta_q^{q-3} + \ldots +1)e_1+\ldots+(\zeta_q+1)e_{q-2}+e_{q-1}.
\end{eqnarray*}

Now since $A \cap L$ is a $\langle x \rangle$-submodule of $A$, $\langle x \rangle$ fixes the line through $f$ in the dual space $A^{*}$.
Multiplying $f$ by $\zeta_q-1$ we get the point $$(\zeta_q^{q-1}-1)e_1+(\zeta_q^{q-2}-1)e_2+\ldots+(\zeta_q-1)e_{q-1}$$ and adding the
functional $e_1+\ldots+e_q$ which is zero on $A$ to it, we obtain another point
\begin{eqnarray*}
\tilde{f} &=& \zeta_q^{q-1}e_1+\zeta_q^{q-2}e_2+\ldots+\zeta_qe_{q-1}+e_q \\
          &=& \sum_{i=1}^{q}\zeta_q^{q-i}e_i.
\end{eqnarray*}

Just as $Sym(q)$ acts on the module $B$ by permuting the basis elements $\theta_i$ via a right action, it acts on the dual space $B^*$ via a
left action by permuting the $e_i$'s, and hence also on $A^*=B^*/\mathbb{F}_p(e_1+\ldots+e_q)$. \vspace{12pt}

Suppose $\sigma \in Sym(q)$ leaves the line through $\tilde{f}$ fixed. That is, $\sigma \tilde{f}=\nu \tilde{f}$ for some $\nu \in
\mathbb{F}_p$. But then
$$\sigma \tilde{f} = \sum_{i=1}^{q}\zeta_q^{q-i}e_{\sigma(i)} = \sum_{i=1}^{q}\zeta_q^{q-\sigma^{-1}(i)}e_i=\sum_{i=1}^{q}\nu\zeta_q^{q-i}e_i=\nu\tilde{f}.$$
Therefore $\zeta_q^{q-\sigma^{-1}(i)}=\nu\zeta_q^{q-i}$ for all $1 \leq i \leq q$ and so $\zeta^{i-\sigma^{-1}(i)}=\nu$. That is, the difference
between $\sigma(i)$ and $i$ is constant, so $\sigma$ is a power of the $q$-cycle $(1 \ 2 \ldots \ q)$. Therefore $\langle x \rangle$ cannot fix
the line in $A^*$ through $f$ and hence cannot stabilize $A \cap L$ which is a contradiction. Therefore no core-free subgroup $L$ can exist and
so $\mu(G)=pq$.\vspace{12pt}

Now suppose $q=3$ and $p \equiv 1$ mod $3$ so that the group $G(=G(p,p,3))$ is isomorphic to $(C_p \times C_p) \rtimes Sym(3)$. As before, let
$c_1, c_2$ generate the base group $A$ and $a=(1 \ 2), b=(1 \ 2 \ 3)$ generate $Sym(3)$. Also as before, the action of $a$ and $b$ on the base
group $A$ induce a two dimensional $Sym(3)$-module structure on $A$. Thus
$$c_1^a=c_1^{-1}, \quad c_2^a=c_1c_2, \quad c_1^b=c_2, \quad c_2^b=c_1^{-1}c_2^{-1}.$$ Let $\zeta_3$ be the
primitive cube root of unity in $\mathbb{F}_p$, so that $\zeta_3^2+\zeta_3+1=0$. Consider the element $c_1c_2^{-\zeta_3}$. We have
$$(c_1c_2^{-\zeta_3})^b=c_1^{\zeta_3}c_2^{\zeta_3+1}=(c_1c_2^{-\zeta_3})^{\zeta_3},$$ so $c_1c_2^{-\zeta_3}$ is an eigenvector for $b$ with
eigenvalue $\zeta_3$. It is easily verified that $c_1c_2^{-\zeta_3}$ is not an eigenvector for $a$ and so the subgroup $L=\langle
c_1c_2^{-\zeta_3}, b \rangle$ forms a core-free subgroup of $G$ of order $3p$. Since $G$ has order $6p$, we have $|G:L|=2p$, so $\mu(G)=2p$.
\vspace{12pt}

Combining this with the previous arguments we have proved:

\begin{theorem} \label{theorem:result1}
Let $p$ and $q$ be odd primes with $p>q$. Then
\begin{displaymath}
\mu(G(p,p,q))= \left\{ \begin{array}{ll}
pq &\text{if} \quad q\geq 5, \quad \text{or} \quad q=3 \quad \text{and} \quad p \equiv 2 \ \text{mod} \ 3 \\
2p & \text{if} \quad  q=3 \quad \text{and} \quad p \equiv 1 \ \text{mod} \ 3. \\
\end{array} \right.
\end{displaymath}
\end{theorem}

\subsection{$\mu(G(p,p,q))$ for $p < q$} \label{section:p<q}

In this subsection, we assume $p<q$ and exploit the action of $b$ on $A$ given in Equation (\ref{equation:baction}). Form the group
$$H:=\langle c_1, c_2, \ldots, c_{q-1}, b \rangle= A \rtimes C_q$$ and treat $A$ as a cyclic $\langle b \rangle$-module. To
prove that the minimal degree of $G$ is $pq$ in this case, we prove it for this proper subgroup $H$. \vspace{12pt}

Recall that the minimal polynomial for the action of $b$ on $A$ is the cyclotomic polynomial $Q_q(\lambda)$. Thus the irreducibility of the
cyclic $\langle b \rangle$-module $A$ depends upon the irreducibility of $Q_q(\lambda)$ over $\mathbb{F}_p$. That is there is a one to one correspondence
between the irreducible factors of $Q_{q}(\lambda)$ and the irreducible submodules of $A$. Moreover, the irreducible factors of $A$ are
precisely the minimal normal subgroups of $H$ as the proposition below proves.

\begin{proposition}
Let $A=A_1 \times A_2 \times \ldots \times A_l$ be the decomposition of $A$ into irreducible $C_q$-modules. Then the $A_i$ are the minimal
normal subgroups of $H$ of order $p^d$, where $l=\frac{q-1}{d}$ and $d$ is the multiplicative order of $p$ in $\mathbb{F}_q^*$.
\end{proposition}

\begin{proof}
Certainly since the $A_i$ are irreducible $C_q$-modules, they are minimal normal subgroups of $H$. That the order of these groups is $p^d$
follows directly from Theorem \ref{theorem:irreduciblesubmodules} since every irreducible factor of $Q_q(\lambda)$ is of degree $d$ which
implies that every irreducible submodule of $A$ is of dimension $d$. \vspace{12pt}

Conversely, let $N$ be a non-trivial normal subgroup of $H$. One can easily adapt the proof of Corollary \ref{corollary:minnorsubgroup} to show
that $N$ contains $A_i$ for some $i$. Thus the $A_i$ are indeed the minimal normal subgroups of $G$.
\end{proof}

This proposition allows us to determine the structure of the subgroups of $H$.

\begin{proposition} \label{proposition:subgroups}
Let $K$ be a proper subgroup of $H$ whose order is divisible by $q$. Then $K \cong C_q$ or $K=\core(K) \rtimes \langle x \rangle$ where $x$ is
an element of order $q$.
\end{proposition}

\begin{proof}
Suppose $K$ is not isomorphic to $C_q$. Therefore $|K|=p^kq$ for $1 \leq k \leq q-1$ and so $K$ contains Sylow $p$-subgroups of order $p^k$.
\vspace{12pt}

By Sylow's Theorem, the number of Sylow $p$-subgroups divides $q$ and is congruent to $1$ mod $p$. Suppose there are $q$ Sylow $p$-subgroups and
let $T_1$ and $T_2$ be two distinct such subgroups. Then they are both $p$-subgroups of $H$ and are thus contained in $A$ since $A$ is the
unique Sylow $p$-subgroup of $H$. However this implies $\langle T_1,T_2 \rangle$ is a $p$-subgroup of $H$ contained in $K$, which properly
contains both $T_1$ and $T_2$, contradicting that they are maximal $p$-subgroups of $K$. \vspace{12pt}

Therefore there is a unique Sylow $p$-subgroup of $K$ and we may write $$K=Syl_{p}(K) \rtimes \langle x \rangle,$$ where $x$ is an element of
order $q$. Observe that $x=ab^j$ where $a \in A$ and $1 \leq j \leq q-1$. Since $Syl_{p}(K)$ is normal in $K$ we have
$$Syl_{p}(K)=Syl_{p}(K)^x=Syl_{p}(K)^{ab^j}=Syl_{p}(K)^{b^j},$$ which shows that $Syl_{p}(K)$ is normal in $H$.\vspace{12pt}

Now since $K$ does not contain $A$, there exists a minimal normal subgroup $A_i$ that intersects trivially with $K$. This implies that $K$ is
not normal in $H$ since for all $a_i \in A_i$, $x^{a_i}=a_{i}^{-1}a_{i}^{'}x$ for some $a_{i}^{'} \in A_i$. Since $Syl_{p}(K)$ is normal in $H$
and maximal in $K$, it now follows that $\core(K)=Syl_{p}(K)$ and so $K=\core(K) \rtimes \langle x \rangle$.
\end{proof}

\begin{corollary} \label{corollary:coreextend}
Let $K$ be a subgroup of $H$ whose order is divisible by $q$. Suppose $T$ is a product of minimal normal subgroups $A_i$ which are not contained
in $K$. Then $\core(KT)=\core(K)T$ is an internal direct product, unless $KT=H$ in which case $\core(KT)=H$.
\end{corollary}

\begin{proof}
Assume that $KT \neq H$. We have by Proposition \ref{proposition:subgroups} that $K=\core(K) \rtimes \langle x \rangle$ where $x$ is an element
of order $q$. Since no $A_i$ is contained in $K$, we have $\core(K) \cap T =\{1\}$. Both are subgroups of $A$ and so they pairwise commute
giving that $\core(K)T$ is an internal direct product. \vspace{12pt}

We may write the group $KT$ as $(\core(K) \times T) \rtimes \langle x \rangle$ and so again by Proposition \ref{proposition:subgroups}, we have
$\core(KT)=\core(K)T$.
\end{proof}

These results give us a clear picture as to the normal structure of $H$. We now proceed to show that the smallest intransitive permutation
representation of $H$ must consist entirely of codimension $1$ subgroups.

\begin{lemma} \label{lemma:replacement}
Let $L$ be a subgroup of $H$ whose order is divisible by $q$ and which is of codimension at least $2$. Then there exist two subgroups of $H$,
$L_1$ and $L_2$, which properly contain $L$ such that $$\core(L_1) \cap \core(L_2) = \core(L)$$ and $$\frac{1}{|L_1|}+\frac{1}{|L_2|} <
\frac{1}{|L|}.$$
\end{lemma}

\begin{proof}
Write $L=\core(L) \rtimes \langle x \rangle$ where $x$ is an element of order $q$. Since $L$ is of codimension at least $2$, there are two
distinct minimal normal subgroups $N_1$ and $N_2$ which intersect $L$ trivially.\vspace{12pt}

Define two subgroups $L_1=LN_1$ and $L_2=LN_2$. By Corollary \ref{corollary:coreextend}, we have $\core(L_1)=\core(L)N_1$ and
$\core(L_2)=\core(L)N_2$ and it is clear that $\core(L_1) \cap \core(L_2) = \core(L)$.\vspace{12pt}

Now $L$ is a proper subgroup of both $L_1$ and $L_2$ of index at least $p$, where $p > 2$. Therefore
$$\frac{1}{|L_1|}+\frac{1}{|L_2|} < (\frac{1}{2}+\frac{1}{2})\frac{1}{|L|}=\frac{1}{|L|}.$$ In fact equality
can only occur if $p=2$ and $H$ contains more than one central involution.
\end{proof}

Notice that the groups $L_1$ and $L_2$ above have dimension $\dim_{G}(L)+1$.

\begin{theorem}\label{theorem:codimension1}
The faithful intransitive permutation representation of $H$ of smallest degree is given by a collection of codimension $1$ subgroups.
\end{theorem}

\begin{proof}
Let $\mathscr{R}=\{L_1,L_2,\ldots, L_n \}$ be a faithful collection of subgroups $H$ where $\core(L_i)$ is non-trivial for all $i$. If there is
an $i$ for which $\codim_{G}(L_i)=0$, then $\core(L_i)=\soc(H)$ and so
$$\{1\}=\bigcap_{j \neq i}^{n}\core(L_j) \cap \core(L_i)$$ implies $\mathscr{R}\setminus{L_i}$ is also a
faithful representation. \vspace{12pt}

Otherwise we suppose there is an $L_i$ which has codimension at least $2$. Then as is Lemma \ref{lemma:replacement} we may replace it with two
subgroups such that the new collection remains faithful and of smaller degree. We repeat this process until the collection contains only
codimension $1$ subgroups of $H$ and this collection will be of least degree.
\end{proof}

We are now in a position to calculate the minimal degree of $H$. We wish to show that every minimal faithful permutation representation is
necessarily transitive, so we must rule out all the intransitive representations. By Theorem \ref{theorem:codimension1}, we only need to rule
the possibility that there is a minimal faithful intransitive representation consisting of codimension $1$ subgroups. We first observe the
following trivial lemma.

\begin{lemma} \label{lemma:numbers}
Let $r$ and $n$ be integers such that $r \geq 3$ and $n \geq 2$. Then $r^{n-1} > n$.
\end{lemma}

\begin{theorem}
Every minimal faithful permutation representation of $H$ is transitive and of degree $pq$.
\end{theorem}

\begin{proof}
We have $H=A \rtimes \langle b \rangle$. If $A$ is an irreducible $C_q$-module, then it follows that $A$ is the unique minimal normal subgroup
of $H$. Therefore in this case, it is clear that every minimal faithful representation of $H$ is transitive. \vspace{12pt}

So suppose that $A$ is a direct product of irreducible $C_q$-modules and write $$H= (\prod_{i=1}^{l}A_i) \rtimes \langle b \rangle,$$ where
$l=\frac{q-1}{d}$. Define subgroups $L_j$ for each $j \in \{1,2,\ldots,l\}$ by
$$L_j= (\prod_{i\neq j}^{l}A_i) \rtimes \langle b \rangle.$$ Then each $L_j$ is a codimension $1$ subgroup of
$H$ with index
$$|H:L_j|=\frac{p^{q-1}}{p^{d(\frac{q-1}{d}-1)}q}=\frac{p^{q-1}}{p^{q-1-d}}=p^d.$$ We have by Proposition \ref{proposition:subgroups} that $\core(L_j)=\prod_{i\neq
j}^{l}A_i$ and so $$\bigcap_{j=1}^{l}\core(L_j)=\bigcap_{j=1}^{l}(\prod_{i\neq j}^{l}A_i)=\{1\},$$ thus the collection
$\mathscr{R}=\{L_1,L_2,\ldots,L_l\}$ is faithful. \vspace{12pt}

Suppose that $\mathscr{R}$ yields a minimal faithful representation. We have
$$\deg(\mathscr{R})=\sum_{j=1}^{l}p^d=\big(\frac{q-1}{d}\big)p^d$$ and since $p(q-1) \leq \mu(H) \leq pq$, this
gives $$p(q-1) \leq \big(\frac{q-1}{d}\big)p^d \leq pq.$$ This simplifies to $(q-1) \leq (\frac{q-1}{d})p^{d-1}\leq q$ and since $q$ is prime,
we are forced to have $q-1 = (\frac{q-1}{d})p^{d-1}$. Therefore $p^{d-1}=d,$ but $d \geq 2$ and $p \geq 3$, so this contradicts Lemma
\ref{lemma:numbers}.\vspace{12pt}

So we have have ruled out the possibility that a minimal faithful representation of $H$ can be given by codimension $1$ subgroups and so Theorem
\ref{theorem:codimension1} implies that every minimal faithful representation must be transitive and thus be given by a core-free
subgroup.\vspace{12pt}

Let $L$ be a core-free subgroup of $H$. We claim that $|H:L| \geq pq$. Suppose for a contradiction that $\core(L)=\{1\}$ and that $|H:L| <pq$.
Since $|H|=p^{q-1}q$, we must have that $|L| > p^{q-2}$. Therefore $|L|=p^{q-1}$ or $q$ divides $|L|$.\vspace{12pt}

The case $|L|=p^{q-1}$ can be ruled out immediately since this implies that $L=A$ which is normal in $H$, contradicting that $L$ is core-free.
So suppose $q$ divides $|L|$. \vspace{12pt}

By Proposition \ref{proposition:subgroups}, we have that $L= \core(L) \rtimes \langle x \rangle$ where $x$ has order $q$. Since we are assuming
that $\core(L)$ is trivial, we must have $L=\langle x \rangle$. Therefore the index of $L$ in $H$ is $p^{q-1}$ and by assumption, $|H:L|=p^{q-1}
<pq$. On the other hand, $p(q-1) \leq \mu(H)$, so $ p(q-1)\leq p^{q-1} < pq$ and since $q$ is prime, this gives $ (q-1)= p^{q-2}$. Since $p$ and
$q$ are odd primes with $p < q$ we have $p$ is at least $3$, $q$ is at least $5$ so this contradicts Lemma \ref{lemma:numbers}. Therefore $L$
has index at least $pq$ as claimed. \vspace{12pt}

So we have shown that any faithful representation of $H$ is transitive of degree at least $pq$ and since $\mu(H) \leq pq$, we have proved
$\mu(H)=pq$.
\end{proof}

It is now clear that we have proved
\begin{theorem} \label{theorem:result2}
If $p$ and $q$ are odd primes with $p<q$, then $\mu(G(p,p,q))=pq$.
\end{theorem}

\subsection{$\mu(G(p,p,p))$}

Throughout this section, $p$ will still denote an odd prime. Our method for calculating $\mu(G(p,p,p))$ is similar to the previous case in that
we calculate the minimal degree of $H=\langle c_1,c_2,\ldots,c_{p-1},b \rangle$ where the $c_i$ generate the group $A(p,p,p)$ and $b=(1 \ 2
\ldots \ p) \in Sym(p)$. However, the structure of the group $H$ differs in this case.

\begin{proposition}
The centre $Z(H)$ of $H$ is $\langle c_1c_2^2 \ldots c_k^k \ldots c_{p-1}^{p-1} \rangle \cong C_p$.
\end{proposition}

\begin{proof}
Suppose $\alpha=c_1^{i_1}c_2^{i_2}\ldots c_{p-1}^{i_{p-1}}b^k$ lies in $Z(H)$ for some $i_1,\ldots i_{p-1},k \in \mathbb{F}_p$. Then $\alpha$
commutes with $c_1$ say, and thus $c_1$ commutes with $b^k$ forcing $k=0$. We also have that $\alpha$ commutes with $b$ and so $\alpha^b=\alpha$
where
$$\alpha^b=c_1^{-i_{p-1}}c_2^{i_1-i_{p-1}}\ldots c_k^{i_k-i_{p-1}}\ldots c_{p-1}^{i_{p-2}-i_{p-1}}.$$ Equating exponents gives $i_s=si_1$.
Without loss of generality we may set $i_1=1$ this giving $\alpha= c_1c_2^2 \ldots c_k^k \ldots c_{p-1}^{p-1}$. Any other choice of $i_1$ would
simply yield a power of $\alpha$. Therefore $Z(H)=\langle c_1c_2^2 \ldots c_k^k \ldots c_{p-1}^{p-1} \rangle \cong C_p$.
\end{proof}

Now $H$ is a $p$-group so any non-trivial normal subgroup of $H$ intersects and hence contains the centre since $Z(H)$ is a copy of $C_p$. This
immediately forces any minimal faithful representation of $H$ to be transitive as before. We may now prove

\begin{proposition} \label{proposition:result3}
$\mu(H)=\mu(G(p,p,p))=p^2$.
\end{proposition}

\begin{proof}
We only need to show that $\mu(H)=p^2$. We show that any core-free subgroup $L$ of $H$ has index at least $p^2$. Suppose this were false. Then
since $|H|=p^p$, $|L| \geq p^{p-1}$. However this gives $|H:L| \leq p$ so $\mu(H) \leq p$ contradicting that $p(p-1)=\mu(A(p,p,p)) \leq \mu(H)$.
\vspace{12pt}

Again $H$ is a proper subgroup of $G(p,p,p)$ which is a proper subgroup of $C_p \wr Sym(p)$. Hence $$p^2 \leq \mu(H) \leq \mu(G(p,p,p)) \leq
p^2,$$ proving the proposition.
\end{proof}

\section{Minimal Degrees of Direct Products} \label{section:directproduct}
 One of the themes of Johnson and Wright's work was to establish conditions for when equality in (\ref{eq:directsum}) holds. The next result is due to Wright \cite{W75}.

\begin{theorem} \label{theorem:nilpotent}
Let $G$ and $H$ be non-trivial nilpotent groups. Then $\mu(G \times H)=\mu(G) + \mu(H)$.
\end{theorem}

Further in \cite{W75}, Wright constructed a class of finite groups $\mathscr{C}$ with the property that for all $G \in \mathscr{C}$, there
exists a nilpotent subgroup $G_1$ of $G$ such that $\mu(G_1)=\mu(G)$. It is a consequence of Theorem (\ref{theorem:nilpotent}) that
$\mathscr{C}$ is closed under direct products and so equality in (\ref{eq:directsum}) holds for any two groups $H,K \in \mathscr{C}$. Wright
proved that $\mathscr{C}$ contains all nilpotent, symmetric, alternating and dihedral groups, however the extent of it is still an open problem.
In \cite{EP88}, Easdown and Praeger showed that equality in (\ref{eq:directsum}) holds for all finite simple groups. \vspace{12pt}

In the closing remarks of \cite{W75}, Wright flagged the question whether equality in (\ref{eq:directsum}) holds for all finite groups. The
referee to that paper then provided an example where strict inequality holds. This formed the basis for considering the minimal degree of these
complex reflection groups and we now demonstrate that strict inequality in (\ref{eq:directsum}) holds when the groups $G(p,p,q)$ are taken in a
direct product. \vspace{12pt}

Let $p$ and $q$ be distinct odd primes as above and let $W=C_p \wr Sym(q)$ be the standard wreath product of the cyclic group of order $p$ by
the symmetric group of degree $q$. Let $\gamma_1, \gamma_2, \ldots, \gamma_q$ be generators for the base group of $W$ and let $a=(1 \ 2)$ and
$b=(1 \ 2 \ldots \ q)$ be the generators for $Sym(q)$ acting coordinate wise on the base group. Let $\gamma=\gamma_1\gamma_2\ldots\gamma_q$,
then $\gamma$ is contained in the centre of $W$. \vspace{12pt}

Let $c_1=\gamma_1\gamma_2^{-1},c_2=\gamma_2\gamma_3^{-1},\ldots, c_{q-1}=\gamma_{q-1}\gamma_{q}^{-1}$. Then it follows that $\langle
c_1,c_2,\ldots,c_{q-1},a,b \rangle$ is isomorphic to $G(p,p,q)$. \vspace{12pt}

Now we claim that $\langle \gamma \rangle \cap G(p,p,q)= \{1\}$. For if this intersection were non-trivial, then we could write
$$\gamma=c_1^{i_1}c_2^{i_2}\ldots c_{q-1}^{i_{q-1}}$$ for some $i_1,i_2,\ldots,i_{q-1}$ in $\mathbb{F}_p$. But upon expressing each side as a
product of the $\gamma_i$'s we have
$$\gamma_1\gamma_2\ldots\gamma_q=\gamma_1^{i_1}\gamma_2^{i_2-i_1}\ldots\gamma_{q-1}^{i_{q-1}-i_{q-2}}\gamma_q^{-i_{q-1}},$$ and equating
exponents gives a contradiction. \vspace{12pt}

It now follows that $W$ is the internal direct product of $\langle \gamma \rangle$ and $G(p,p,q)$ and that $\mu(\langle \gamma \rangle \times
G(p,p,q))=\mu(G(p,p,q))=pq$. Therefore $$ pq=\mu(\langle \gamma \rangle \times G(p,p,q)) <  \mu(\langle \gamma \rangle)+\mu(G(p,p,q))=p+pq$$ and
we have a class of groups for when strict inequality in (\ref{eq:directsum}) holds.\vspace{12pt}

In the case where $p=q$, the centralizer of $G(p,p,p)$ in $Sym(p)$ is properly contained in $G(p,p,p)$. So (\ref{eq:directsum}) is an equality
in this case.

\section{Acknowledgements}
The author acknowledges ideas from a preprint (currently being revised) by Ben Elias, Lior Silberman and Ramin Takloo-Bighash which gave
assistance to results in Section \ref{section:p<q}.\vspace{12pt}

The author would also like to thank his supervisor David Easdown and associate supervisor Anthony Henderson for many helpful discussions and
their expert guidance.


\bibliographystyle{plain}

\vspace{12pt}

\end{document}